\documentclass[11pt]{amsart}
\usepackage{amscd,graphicx}

\newtheorem{theo}{Theorem}
\newtheorem{lemm}[theo]{Lemma}
\newtheorem{step}{Step}

\theoremstyle{remark}
\newtheorem*{remark}{Remark}

\newcommand{\Z}{{\mathbb Z}}
\newcommand{\R}{{\mathbb R}}
\newcommand{\C}{{\mathbb C}}
\newcommand{\N}{{\mathfrak N}}
\newcommand{\w}{{\overline w}}
\newcommand{\pa}{{\partial}}
\newcommand{\rank}{\operatorname{rank}}
\newcommand{\krn}{\operatorname{ker}}
\newcommand{\img}{\operatorname{im}}
\renewcommand{\mod}{\operatorname{mod}}
\newcommand{\Hom}{\operatorname{Hom}}
\newcommand{\Ext}{\operatorname{Ext}}

\begin{document}

\title[Cobordisms of fold maps\dots]
{Cobordisms of fold maps and maps with prescribed number of cusps}

\author{Tobias Ekholm}
\email{tekholm{\@@}usc.edu}
\address{USC Department of mathematics, 3620 S Vermont Ave, Los Angeles, CA 90089}
\author{Andr\'as Sz\H ucs}
\email{szucs{\@@}cs.elte.hu}
\address{Department of Analysis, E{\"o}tv{\"o}s Lor{\'a}nd University (ELTE),
1088 Budapest, HUNGARY}
\author{Tam\'as Terpai}
\email{terpai{\@@}cs.elte.hu}
\address{Department of Analysis, E{\"o}tv{\"o}s Lor{\'a}nd University (ELTE),
1088 Budapest, HUNGARY}

\date{}
\thanks{TE is an Alfred P. Sloan Research Fellow and acknowledges support from NSF-grant DMS-0505076.}

\begin{abstract}
A generic smooth map of a closed $2k$-manifold into $(3k-1)$-space has a finite number of cusps ($\Sigma^{1,1}$-singularities). We determine the possible numbers of cusps of such maps. A fold map is a map with singular set consisting of only fold singularities ($\Sigma^{1,0}$-singularities). Two fold maps are fold bordant if there are cobordisms between their source- and target manifolds with a fold map extending the two maps between the boundaries, if the two targets agree and the target cobordism can be taken as a product with a unit interval then the maps are fold cobordant. We compute the cobordism groups of fold maps of $(2k-1)$-manifolds into $(3k-2)$-space.  Analogous cobordism semi-groups for arbitrary closed $(3k-2)$-dimensional target manifolds are endowed with Abelian group structures and described. Fold bordism groups in the same dimensions are described as well.
\end{abstract}

\maketitle

\section{Introduction}\label{sec:1}
Let $M$ be a manifold, $\dim(M)=n$, and let $f\colon M\to\R^{n+k}$ be a smooth map. A point $p\in M$ is a {\em singular point} of $f$ if the rank of the differential of $f$ at $p$, $\rank(df_p)$, is smaller than $n$. We let $\Sigma(f)\subset M$ denote the set of singular points of $f$. If $f$ is a stable map (i.e., if the orbit of $f$ under left-right action by diffeomorphisms is open in the space of smooth maps $M\to\R^{n+k}$) then $\Sigma(f)$ is naturally stratified
$$
\Sigma(f)=\Sigma^1(f)\cup\Sigma^2(f)\cup\dots,
$$
where $\Sigma^j(f)$ is the set of singular points $p$ with $\rank(df_p)=n-j$. Studying the restriction of $f$ to the strata one obtains further stratifications of $\Sigma^j(f)$ called the Thom-Boardman stratification, see \cite{Bo}. In this paper we will be concerned only with the simplest singularities in the Thom-Boardman classification, we give a brief description of these.

Consider a stable map $f\colon M\to\R^{n+k}$ with $\Sigma(f)=\Sigma^1(f)$.
The $1$-dimensional kernel of the differential $\krn(df)$ is a line field in
$T_{\Sigma(f)} M$ (i.e., in the restriction of the tangent bundle of $M$ to $\Sigma(f)$). Let $p\in\Sigma(f)$. If $\krn(df)$ is not tangent to $\Sigma(f)$ at $p$ then $p$ is a {\em $\Sigma^{1,0}$-point}. If $\krn(df)$ is tangent to $\Sigma(f)$ at $p$ then $p$ is a {\em
$\Sigma^{1,1}$-point}. A $\Sigma^{1,1}$-point $p$ is a {\em cusp} if the tangency of $\krn(df)$ with $\Sigma(f)$ is transverse at $p$.

Let $\Sigma^{1,0}(f)$ and $\Sigma^{1,1}(f)$ denote the sets of $\Sigma^{1,0}$-points and $\Sigma^{1,1}$-points, respectively. It is a straightforward consequence of the jet transversality theorem that $\Sigma^{1,0}(f)\subset M$ has codimension $k+1$, that $\Sigma^{1,1}(f)\subset M$ has codimension $2(k+1)$, and that the set of points where $\krn(df)$ has a non-transverse tangency with $\Sigma(f)$ (the set of $\Sigma^{1,1}(f)$-points which are not cusps) has codimension $3(k+1)$.
We say that a stable map $f\colon M\to\R^{n+k}$ is a {\em fold map} if $\Sigma(f)=\Sigma^{1,0}(f)$.

Let $M_0$ and $M_1$ be closed $n$-manifolds. Two fold maps $f_j\colon M_j\to\R^{n+k}$, $j=0,1$, are {\em fold cobordant} if there exists a cobordism $W$ between $M_0$ and $M_1$ (i.e., $W$ is an $(n+1)$-manifold with $\pa W=M_0\sqcup M_1$) and a fold map $F\colon W\to\R^{n+k}\times[0,1]$ such that $F|M_j=f_j\times\{j\}$, $j=0,1$. There are analogous notions for oriented manifolds: two fold maps $f_j\colon M_j\to\R^{n+k}$, $j=0,1$, of closed oriented manifolds are {\em oriented fold cobordant} if there exists an oriented cobordism $W$ (i.e., $W$ is an oriented $(n+1)$-manifold with $\pa W=-M_0\cup M_1$) and a fold map $F\colon W\to\R^{n+k}\times[0,1]$ extending $f_j\times\{j\}$, $j=0,1$, on $\pa W$. Following \cite{A} and \cite{R-Sz} we define {\em the cobordism group of (oriented) fold maps of $n$-manifolds into $\R^{n+k}$}, $\Sigma^{1,0}(n,k)$ ($\Sigma_{\rm so}^{1,0}(n,k)$) to be the set of fold maps of (oriented) $n$-manifolds up to (oriented) fold cobordism, with group operation induced by disjoint union.

It follows immediately from the computation of the codimension of $\Sigma^{1,1}(f)$ above that if $2k+1>n$ then the forgetful morphisms from $\Sigma^{1,0}(n,k)$ and $\Sigma^{1,0}_{\rm so}(n,k)$ which forget the map are isomorphisms to the cobordism group of $n$-manifolds $\mathfrak N_n$ and to the oriented cobordism group of $n$-manifolds $\Omega_n$, respectively. We compute the cobordism groups of fold maps for the largest codimension where transversality arguments do not imply that
the forgetful morphisms are isomorphisms i.e., for $2k+1=n$.

\begin{theo}\label{t1}
The cobordism groups of fold maps are as follows.
\begin{itemize}
\item[{\rm (a)}] For any $k\ge 0$,
$$
\Sigma^{1,0}(2k+1, k) \approx \N_{2k+1},
$$
\item[{$\rm (b_0)$}]
$$
\Sigma^{1,0}_{\rm so}(1,0)\approx\Z_2
$$
\item[{$\rm (b_1)$}]
$$
\Sigma_{\rm so}^{1,0}(5,2) \approx \Omega_5\oplus\Z_2\approx\Z_2\oplus\Z_2
$$
\item[{$\rm (b_\ast)$}] For any $m\ge 2$
$$
\Sigma^{1,0}_{\rm so}(4m+1,2m)\approx
\Omega_{4m+1}.
$$
\item[{\rm (c)}] For any $m\ge 1$,
$$
\Sigma^{1,0}_{\rm so}(4m-1,2m-1)\approx\Omega_{4m-1}\oplus\Z_{3^t},
$$
where $t = \min\{j\mid \alpha_3(2m + j)\leq 3j\}$
and $\alpha_3(x)$ denotes the sum of digits of the integer $x$
in triadic system. For example, $\Sigma^{1,0}_{\rm so}(3,1)\approx\Z_3$ and
$\Sigma^{1,0}_{\rm so}(33, 15)\approx\Omega_{33}\oplus\Z_9$
\end{itemize}
\end{theo}

Theorem \ref{t1} is proved in Subsection \ref{s:bgfm}. For $n$ odd, the group $\Omega_n$ is isomorphic to $\Z_2\oplus\dots\oplus\Z_2$. The number of summands was determined by Wall \cite{W}. The summand $\Z_{3^t}$ in (c) appears as a consequence of a result of Stong \cite{St}, which describes the possible values of the top normal Pontryagin class of an oriented $4m$-manifolds.

Theorem \ref{t1} is a consequence of the following result which determines the possible number of cusps of a stable map $F\colon W\to\R^{3k-1}$ of a closed (oriented) $2k$-manifold $W$. In order to state it we note that the cusps of a stable map of an oriented manifold $W$ of dimension $2k=4m$ are oriented $0$-manifolds, i.e., points with signs. See \cite{Sz2} or \cite{E-Sz}, Section 6, for an alternative  approach. In this case we say that the {\em algebraic number of cusps} is the sum of signs over the cusps.

\begin{theo}\label{l2}
${ }$
\begin{itemize}
\item[{\rm (a)}] For any integers $c\ge 0$ and $k>0$, there exist a (possibly non-orientable)
closed manifold $W$ of dimension $2k$ and a stable map $F\colon W\to \R^{3k-1}$ with $c$ cusps.
\item[{$\rm (b_{01})$}] If $W$ is an orientable closed $(4m+2)$-manifold, $m=0,1$
and $F\colon W\to\R^{6m+2}$ is a stable map then the number of cusps of $F$ is even and for any even integer $c\ge 0$ there exists a map $F'\colon W\to W$ with $c$ cusps.
\item[{$\rm (b_\ast)$}] For any integers $c\ge 0$ and $m\ge 2$, there exist an orientable closed manifold $W$ of dimension $4m+2$ and a stable map $F\colon W \to \R^{6m+2}$ with $c$  cusps.
\item[{\rm (c)}] If $W$ is a closed oriented manifold of dimension $4m$, $m>0$, and $F\colon W\to\R^{6m-1}$ is a stable map then the algebraic number of cusps of $F$ equals the normal Pontryagin number ${\overline p}_m[W]$ of $W$. Moreover for any $c$, such that $c= |{\overline p}_m[W]|+2r$ for some integer $r\ge 0$, there exists a stable map $F'\colon W\to\R^{6m-1}$ with $c$ cusps.
\end{itemize}
\end{theo}
Theorem \ref{l2} is proved in Subsection \ref{s:pnc}.

In Sections \ref{sec:bor-cob} and \ref{sec:4}, we consider problems similar to those discussed above replacing the target spaces $\R^{3k-2}$ by closed $(3k-2)$-manifolds. More specifically, in Section
\ref{sec:bor-cob} we allow also the target to change by cobordism: Two stable maps $f_j\colon V_j\to X_j$, $j=0,1$, where $V_j$ are closed $(2k-1)$-manifolds and where $X_j$ are closed $(3k-2)$-manifolds are {\em fold bordant} if there exists cobordisms $W$ and $Y$ with $\pa W= V_0\sqcup V_1$ and $\pa Y = X_0\sqcup X_1$, respectively, and a fold map $F\colon W\to Y$ extending the maps $f_j$, $j=0,1$, on the boundary. Equivalence classes of fold bordant maps naturally form groups which are presented in Theorem \ref{t:last}.

In Section \ref{sec:4}, we study fold cobordism classes of maps of $(2k-1)$-manifolds into arbitrary (but fixed) closed $(3k-2)$-manifolds. In this situation the equivalence classes form a semi-group rather than a group. (There is no natural geometric construction of an inverse.) However, in Theorem \ref{t:lastlast} we show how to endow these semi-groups with Abelian group structures and describe the corresponding groups.

\section{Cusp cancelation}\label{sec:3/2}
In this section we show how to eliminate pairs of cusps. According to Morin
\cite{M} if $F\colon W\to\R^{3k+2}$ is a stable map of a $(2k+2)$-manifold with a cusp at $p\in W$ then there exist coordinates
$$
(t,s,x)=(t_1,\dots,t_{2k},s,x)\in\R^{2k+2},\quad \text{around }p\in W
$$
(i.e., $p$ corresponds to the origin) and
$$
(y,z,u)=(y_1,\dots,y_{2k+1},z_1,\dots,z_k,u)\in\R^{3k+2},\quad\text{around }F(p)\in\R^{3k+2}
$$
such that $(y,z,u)=F(t,s,x)$ is given by the expression
\begin{align*}
&y_j=t_j,\quad j=1,\dots,2k,\\
&y_{2k+1}=s,\\
&z_j=xt_{2j-1}+ x^2 t_{2j},\quad  j=1,\dots,k,\\
&u=x(x^2+s).
\end{align*}

\begin{lemm}\label{l3}
Let $F\colon W\to\R^{3k+2}$, $k>0$, be a stable map of a $(2k+2)$-manifold and
let $p$ and $p'$ be cusps of $F$. If $W$ is oriented and $k$ is odd assume
that the signs of $p$ and $p'$ are different. Assume that there exists an arc
$\alpha$ connecting $p$ to $p'$ in $W$ such that
$\Sigma(F)\cap\alpha=\{p,p'\}$. Then there exists a stable map $F'\colon
W\to\R^{3k+2}$ which agrees with  $F$ outside an arbitrarily small neighborhood of $\alpha$, and such that $\Sigma^{1,1}(F')=\Sigma^{1,1}(F)-\{p,p'\}$.
\end{lemm}

\begin{remark}
For the counterpart of Lemma \ref{l3} in the case $k=0$, see e.g. \cite{Ka}.
\end{remark}

\begin{proof}
Consider the disk
$$
\bigl\{(t,s,x)\colon |t|\le\epsilon, -\epsilon\le s\le 1+\epsilon, |x|\le 1\bigr\}
$$
in coordinates as above and two maps $\phi_1$ and $\phi_2$ from $\R^{2k+2}$ to
$\R^{3k+2}$ which have the same $(y,z)$-components as the map above and for which the $u$-coordinates have the form
$$
(\phi_j)_u(x,s)=x(x^2+\psi_j(s)),\quad j=1,2,
$$
for $|x|<\frac{15}{16}$, where the functions $\psi_j(s)$ satisfy
\begin{align*}
\psi_1(s)&=-s(s-1),\\
\psi_2(s)&= C(s)-s(s-1).
\end{align*}
Here the function $C(s)$ equals $0$ for $-\epsilon\le s\le-\frac34\epsilon$ and for $1+\frac34\epsilon\le s<1+\epsilon$, equals  $-\frac13$ for $-\frac14\epsilon\le s\le 1+\frac14\epsilon$, is decreasing for $-\frac34\epsilon<s<-\frac14\epsilon$, and increasing for $1+\frac14\epsilon<s<1+\frac34\epsilon$. Then the function $x\mapsto(\phi_1)_u(x,s)$ has two critical points for all $s$, and the function $x\mapsto(\phi_2)_u(x,s)$ has two critical points for $-\epsilon\le s<0$ which cancel at $s=0$ and two critical points which are born at $s=1$ and persists for $1< s\le 1+\epsilon$. Interpolating between the two formulas in the region $\frac{15}{16}<|x|<1$ it is easy to arrange that $\phi_1$ and $\phi_2$ agree on the boundary of the disk. Note that $\phi_1$ does not have any cusps whereas $\phi_2$ has two cusps (of opposite signs if $k$ is odd).

To prove the lemma we note that the restriction of $F$ to a neighborhood of $\alpha$ is equivalent to the map $\phi_1$ (under left-right action of diffeomorphisms), hence we can replace $F$ on this neighborhood of $\alpha$ with a map equivalent to $\phi_2$ thereby removing the cusps of $F$ as desired.
\end{proof}

\section{Proofs of Theorems \ref{t1} and \ref{l2}}\label{sec:2}
In this section we first prove Theorem \ref{l2} and then Theorem \ref{t1}.

\subsection{Maps with prescribed number of cusps}\label{s:pnc}
\begin{proof}[Proof of Theorem \ref{l2} (c)]
This was proved in \cite{Sz2}, see also \cite{E-Sz}.
\end{proof}

\begin{proof}[Proof of Theorem \ref{l2} (a)]
The case $k=1$ is well known see e.g. \cite{Ka}. Let $k>1$. If
$F\colon W\to\R^{3k-1}$ is a stable map with at least two cusps then by adding
$1$-handles we can make $W$ connected. Since the codimension of the singular
set is $k+1>1$ we may use Lemma \ref{l3} to remove pairs of cusps.  Thus it
suffices to show that there exists a manifold $W$ and a stable map $F\colon
W\to\R^{3k-1}$ with an odd number of cusps. Let $W=(\R P^2)^{k+1}$ where $\R
P^2$ is the real projective plane. Let $b\colon\R P^2\to\R^3$ be an immersion
(for example the Boy surface), and let $g\colon W\to \R^{3k+3}$ be a
self-transverse immersion regularly homotopic to the product $b\times
{\cdots}\times b$, ($(k+1)$-factors). Let $F: W\to \R^{3k + 2}$ be a composition of $g$ with a generic hyperplane projection. It is well known that $g$ has an odd number of triple points,
see for example \cite{H}. In \cite{Sz2} it was shown that the number of cusps of $F$ is
congruent $\mod 2$ to the number of triple points of $g$.
\end{proof}

\begin{proof}[Proof of Theorem \ref{l2} (${\it b_{01}}$)]
For $m=0$ it is well known that the number of cusps of a stable map $F$ of an
orientable surface $W$ into the plane has an even number of cusps, see e.g.
Theorem 9 in \cite{T}.

For $m=1$, note that $\Omega_6=0$. Therefore, if $F\colon W\to\R^8$ is any stable map of an oriented $6$-manifold, there exists a compact orientable $7$-manifold $N$ with $\pa N=W$ and a stable map $G\colon N\to\R^8\times\R_+$ extending $F$. Now $\Sigma^{1,1}(G)$ is a $1$-manifold with boundary $\Sigma^{1,1}(F)$. This proves that the number of cusps of $F$ is even.
\end{proof}

\begin{proof}[Proof of Theorem \ref{l2} (${\it b_\ast}$)]
As in the proof of (a), we note that by adding a $1$-handle in an orientation preserving manner we can cancel pairs of cusps. Thus it suffices to produce an oriented $(4m+2)$-manifold $W$ and a map $F\colon W\to\R^{6m+2}$, $m\ge 2$, with an odd number of cusps.
To this end let $Y$ be the Dold manifold $(\C P^2\times S^1)/\Z_2$, where $\Z_2$ acts by complex conjugation on $\C P^2$ and by multiplication by $-1$ on $S^1$. This
manifold is orientable and generates $\Omega_5\approx\Z_2$, see \cite{MSt}. Since the manifold $\left(Y\times(\mathbb{R}P^2)^{k-2}\right)^2$ is a square there exists an orientable
manifold cobordant to it, see \cite{W} or \cite{CF}. Let $W$ be such a manifold and let
$F\colon W \to \R^{6k+2}$ be a stable map. We compute the parity of the number
of cusps of $F$ (which we denote by $\#_2\Sigma^{1,1}(F)\in\{0,1\}\approx\Z_2$) using the formula
\begin{equation}\label{eq:nofcusp}
\#_2\Sigma^{1,1}(F)\equiv\left({\w}_{2k+1}^2+{\w}_{2k+2} {\w}_{2k}\right)[W],
\end{equation}
see \cite{Ro2}, or \cite{Ro} in combination with \cite{BH}: \cite{Ro}
contains the complex analogue of this statement and Theorem 6.2 in \cite{BH} shows that the complex version implies the real.
Let $N=Y\times(\R P^2)^{k-2}$. Applying the product formula for normal Stiefel-Whitney classes, we
obtain
\begin{align*}
\w_{2k+1}^2[N\times N]&=\sum_{i_1+j_1=2k+1}\,\,\sum_{i_2+j_2=2k+1}
(\w_{i_1}\times\w_{j_1})(\w_{i_2}\times\w_{j_2})[N\times N]\\
&=\sum_{i_1+j_1=2k+1}\,\,\sum_{i_2+j_2=2k+1}
(\w_{i_1}\,\w_{i_2})[N](\w_{j_1}\,\w_{j_2})[N]\\
&=2\sum_{0\le i\le k}(\w_i\,\w_{2k+1-i})[N](\w_{2k+1-i}\,\w_i)[N]=0.
\end{align*}
Similarly,
\begin{align*}
\w_{2k+2}\w_{2k}[N\times N]
&=\sum_{i_1+j_1=2k+2}\,\,\sum_{i_2+j_2=2k}
(\w_{i_1}\times\w_{j_1})(\w_{i_2}\times\w_{j_2})[N\times N]\\
&=\sum_{i_1+j_1=2k+2}\,\,\sum_{i_2+j_2=2k}
(\w_{i_1}\,\w_{i_2})[N](\w_{j_1}\,\w_{j_2})[N]\\
&=\sum_{0\le i\le 2k}
(\w_{2k+1-i}\,\w_{i})[N](\w_{i+1}\,\w_{2k-i})[N]\\
&=\bigl(\w_{k+1}\,\w_k[N]\bigr)^2=\w_{k+1}\,\w_k[N],
\end{align*}
and, since $\w_i\,\w_j[\mathbb{R}P^2]=1$ if and only if $i=j=1$, we find
\begin{align*}
&\w_{k+1}\,\w_k[N]=\\
&=\sum_{\substack{i_0+\dots+i_{k-2}=k+1\\j_0+\dots+j_{k-2}=k}}
(\w_{i_0}\times\dots\times\w_{i_{k-2}})(\w_{j_0}\times\dots\times\w_{j_{k-2}})
[Y\times(\R P^2)^{k-2}]\\
&=\sum_{\substack{i_0+\dots+i_{k-2}=k+1\\j_0+\dots+j_{k-2}=k}}
\bigl(\w_{i_0}\,\w_{j_0}[Y]\bigr)\bigl(\w_{i_1}\,\w_{j_1}[\R P^2]\bigr)
\dots\bigl(\w_{i_{k-2}}\,\w_{j_{k-2}}[\mathbb{R}P^2]\bigr)\\
&=\w_3\,\w_2[Y]=1.
\end{align*}
To see that the last equality holds we argue as follows. Since $Y$ is
odd-dimensional its top Steifel-Whitney class vanishes: $w_5[Y]=0$. Since $Y$
is orientable $w_1(Y)=0$. A straightforward calculation shows that these two
conditions imply $w_2(Y)=\w_2(Y)$ and $w_3(Y)=\w_3(Y)$. The manifold  $Y$
generates $\Omega_5\approx\N_5$ and thus it is not null-cobordant. Therefore at least one of its Steifel-Whitney numbers must be non-zero and by the two conditions the only possibility is $w_2\,w_3[Y]$. Thus $\w_2\,\w_3[Y]=w_2\,w_3[Y]=1$. This finishes the proof of (b).
\end{proof}

\subsection{Cobordism groups of fold maps}\label{s:bgfm}

\begin{proof}[Proof of Theorem \ref{t1} (a)]
We show that the forgetful morphism
$$
r\colon\Sigma^{1,0}(2k+1,k)\to\N_{2k+1},
$$
which associates to each cobordism class represented by a stable map $f\colon
M\to\R^{3k+1}$ the cobordism class $[M]$ of its source is an isomorphism. The
morphism is surjective by the jet-transversality theorem, so it is sufficient to show that it is
injective. Thus let $f\colon M\to\R^{3k+1}$ be a stable map where $M$ is
null-cobordant. Let $W'$ be a $(2k+2)$-manifold with $\pa W=M$. Extend the map
$f$ to a stable map $F'\colon W'\to\R^{3k+1}\times\R_+$. Assume that $F$ has
$c$ cusps. Using Theorem \ref{l2} (a) we find a closed $(2k+2)$-manifold $W''$
and a stable map $F''\colon W''\to\R^{3k+1}\times\R_+$ which has $c$
cusps. Now adding $1$-handles connecting cusps of the maps $F'$ and $F''$, and
then, applying Lemma \ref{l3}, we produce a
fold  map $F\colon W\to \R^{3k+1}\times \R_+$ where $W$ satisfies $\pa W=M$ and where $F|M=f\times\{0\}$. This shows that $r$ is injective.
\end{proof}

\begin{figure}
\begin{center}
\includegraphics*[width=5cm]{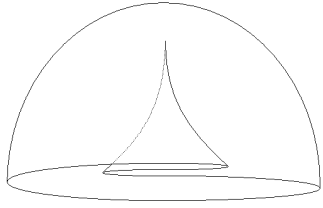}
\caption{A neighborhood of a cusp.}
\label{fig:cusp}
\end{center}
\end{figure}

\begin{proof}[Proof of Theorem \ref{t1} (c)]
If $F\colon W\to\R^{6m-1}$ is a stable map of an oriented $4m$-manifold then let $\#\Sigma^{1,1}(F)$ denote the algebraic number of cusps of $F$. In order to determine the kernel $\krn(r_{\rm so})$ of the forgetful morphism
$$
r_{\rm so}\colon\Sigma^{1,0}_{\rm so}(4m-1,2m-1)\to\Omega_{4m-1},
$$
let $G_m$ be the greatest common divisor of the numbers in the set $\bigl\{\#\Sigma^{1,1}(F)\bigr\}$, where $F$ ranges over all stable maps $F\colon W\to\R^{6k-1}$ and where $W$ ranges over all oriented closed $4m$-manifolds. It follows from Theorem \ref{l2} (c) in combination with a result of Stong \cite{St} that $G_m=3^t$ (with $t$ is as in the formulation of the theorem). Define the map
$$
\Gamma\colon \krn(r_{\rm so})\to\Z_{G_m}
$$
as follows for $\xi\in\krn(r_{\rm so})$. Pick a representative $f\colon
M\to\R^{6m-2}$ of $\xi$. Since $M$ is oriented null-cobordant there exists a compact oriented $4m$-manifold $W$ such that $\pa W=M$. Let $F\colon W\to\R^{6m-2}\times\R_+$ be any stable map which extends $f$ and let
$$
\Gamma(\xi)=\#\Sigma^{1,1}(F)\quad\mod G_m.
$$
In order to see that $\Gamma(\xi)$ is independent of
the choice of  $F$, let $F'\colon W'\to\R^{6m-2}\times\R_+$ be another extension of $f$. If $\hat F$ denotes $F$ composed with the reflection in $\R^{6m-2}$ then $\hat F\cup F'$ gives a map $H$ of the closed oriented manifold $-W\cup_M W'$ into $\R^{6m-1}$ and the algebraic number of cusps of this map is
$$
\#\Sigma^{1,1}(H)=\#\Sigma^{1,1}(F')-\#\Sigma^{1,1}(F)=0\quad\mod G_m.
$$

It is easy to see that $\Gamma(\xi)$ does not depend on the representative $f$
since if $f'$ is another representative then there is a cobordism without cusps connecting them.
Thus $\Gamma(\xi)$ is well defined.

We next show that $\Gamma$ is an isomorphism. First, using a neighborhood of a cusp, see Section \ref{sec:3/2} for formulas and Figure \ref{fig:cusp} (which is also Figure 19 in
 \cite{Ar}) for an illustration, it is easy to construct a stable map of a sphere $S^{4m-1}\to\R^{6m-2}$ which represents an element $\xi$ with $\Gamma(\xi)=1$ showing that $\Gamma$ is surjective. If $f\colon M\to\R^{6m-2}$ is a stable map having an extension $F\colon W\to\R^{6m-2}\times\R_+$ with algebraic number of cusps $c$ divisible by $G_m$ then using Theorem \ref{l2} (c) we find a closed manifold $W'$  and a stable map $F'\colon W'\to\R^{6m-2}\times\R_+$ with algebraic number of cusps equal to $-c$. Using connected sum we may join $F$ and $F'$ to a map of a connected manifold $W''$ with algebraic number of cusps equal to $0$. Since the codimension of the singular set is $2m>0$ there arcs connecting cusp pairs of different signs and Lemma \ref{l3} implies we can cancel all cusps. This shows that $\Gamma$ is injective. We thus have the exact sequence
$$
\begin{CD}
0 @>>> \Z_{G_m} @>>> \Sigma^{1,0}(4m-1,2m-1) @>>> \Omega_{4m-1} @>>>0
\end{CD}
$$
and statement (c) follows since $\Omega_{4m-1}$ is a direct sum of copies of $\Z_2$ and $G_m$ is a power of $3$.
\end{proof}

\begin{proof}[Proof of Theorem \ref{t1} (${\it b_\ast}$)]
Apply the argument in the proof of (a) to the forgetful morphism
$$
r_{\rm so}\colon\Sigma^{1,0}_{\rm so}(4m+1,2m)\to\Omega_{4m+1},
$$
using Theorem \ref{l2} ${\rm (b_\ast)}$ instead of Theorem \ref{l2} (a).
\end{proof}

\begin{proof}[Proof of Theorem \ref{t1} (${\rm b_0}$)]
Using the fact that $\Omega_1=0$, the statement follows from the argument in
the proof of (c) with $G_m=2$. (It is straightforward to give a direct geometric proof of this noting that the fold cobordism class of a map $f\colon S^1\to\R$ is equal to $\frac12(\mu(f)-2)\mod 2$, where $\mu(f)$ is the number of critical points of $f$.)
\end{proof}

\begin{proof}[Proof of Theorem \ref{t1} (${\it b_1}$)]
Applying the argument in the proof of (c) with $G_m=2$ (referring to Theorem \ref{l2} ${\rm (b_{01})}$ instead of Theorem \ref{l2} (c)) we obtain the exact sequence
$$
\begin{CD}
0 @>>> \Z_2 @>>> \Sigma^{1,0}(5,2) @>>> \Omega_{5} @>>>0.
\end{CD}
$$
Noting that $\Omega_5\approx \Z_2$ this leaves the following two possibilities
\begin{itemize}
\item[{(i)}]
$\Sigma^{1,0}_{\rm so}(5,2)\approx\Z_4$, or
\item[{(ii)}]
$\Sigma^{1,0}_{\rm so}(5,2)\approx\Z_2\oplus\Z_2$.
\end{itemize}

If $f\colon M\to\R^7$ is a stable map then let $[f,M]\in\Sigma_{\rm
 so}^{1,0}(5,2)$ denote the fold cobordism class represented by it. If (i) holds
 and if $f\colon M\to\R^7$ is a stable map, where $M$ is not oriented
 null-cobordant, then $[f,M]$ is an element of order $4$. To show that (ii)
 holds it is thus sufficient to establish the existence of a stable map
 $f\colon Y\to\R^7$ such that $Y$ is not oriented null-cobordant and such that
 $[f,Y]$ has order $2$. To this end, as in the proof of Theorem \ref{l2} ${\rm (b_\ast)}$, let $Y$  be the Dold manifold with some fixed orientation. By \cite{Co} there exists an immersion $\tilde f\colon Y\to\R^8$. We show that if $f=\pi\circ\tilde f\colon Y\to\R^7$, where $\pi\colon\R^8\to\R^7$ is a generically chosen projection then $[f,Y]$ is an element of order $2$. The proof has three steps.

\begin{step}\label{step1}
Let $g\colon Y\to\R^7$ be any fold map. Let $H\subset\R^7$ be a hyperplane disjoint from $g(Y)$, let $r\colon\R^7\to\R^7$ be reflection in $H$, and let $-Y$ denote $Y$ with the opposite orientation. Then $[g,Y]+[r\circ g, -Y]=0$.
\end{step}

This follows from the standard construction of a cobordism inverse: if $x_1$ is
a coordinate perpendicular to $H$ and centered on $H$, if $g=(g_1,g')\in\R\times H$, and if $G\colon Y\times[0,\pi]\to\R^8_+$ is given by
$$
G(\theta, y)=\bigl(\sin\theta g_1(y),\cos\theta g_1(y), g'(y)\bigr)\in\R_+\times\R\times H,
$$
then $G$ gives a fold cobordism establishing Step \ref{step1}.

The manifold $Y$ admits an orientation reversing diffeomorphism $A\colon Y\to
 Y$ induced by complex conjugation on the $S^1$-factor thought of as the unit
circle in $\C$. Let $X$ be the mapping torus of $A$, i.e., $X=(Y\times I)/\sim$, where $(y,0)\sim(A(y),1)$. Then $X$ is a non-orientable closed
$6$-manifold. (This is the $6$-dimensional Wall manifold, $X_6$ in the notation of
\cite{W}.) In \cite{W}, the cohomology ring of $X$ and its Stiefel-Whitney classes were computed. This computation in combination with a straightforward calculation imply that the normal Stiefel- Whitney number ${\w}_3^2[X]+ {\w}_2{\w}_4[X]$ equals $0$. Equation \eqref{eq:nofcusp} then implies that the number of cusps of any stable map $F\colon X\to\R^8$ is even.

\begin{step}\label{step2}
Let $g\colon Y\to\R^7$ be any fold map then $[g,Y]+[r\circ g,Y]=0$.
\end{step}

To see this we first define a map $K\colon X\to\R^8$ as follows. Consider $X=(Y\times[0,\frac12]\sqcup Y\times[\frac12,1])/\sim$, where $(y,\frac12)\sim(y,\frac12)$ and $(y,0)\sim(A(y),1)$. Define $K\colon Y\times[0,\frac12]\to\R^8_+$, $K(t,y)=G(2\pi t,y)$ where $G$ is as in Step \ref{step1} and then extend $K$ over $Y\times[\frac12,1]$. Since $K$ does not have any cusps in $Y\times[0,\frac12]$ it follows that $K$ has an even number of cusps in $Y\times[\frac12,1]$. Since $\Omega_5\approx \Z_2$ there exists an orientable $6$-manifold $W$ with $\pa W=Y\sqcup Y$. Let $F\colon W\to\R^8_+$ be a stable map such that the restriction to one boundary component agrees with $g$ and the restriction to the other agrees with $r\circ g$. Joining the maps $F$ and $K|Y\times[\frac12, 1]$ along their common boundary we obtain a map $K\# F$ of an orientable $6$-manifold into $\R^8$. By Theorem \ref{l2} ${\rm (b_{01})}$ this map has an even number of cusps. It follows that $F$ has an even number of cusps and Lemma \ref{l3} then implies that $[g,Y]+[r\circ g,Y]=0$.

\begin{step}\label{step3}
Let $f\colon Y\to\R^7$ be the projection of an immersion $\tilde f\colon Y\to\R^8$ then $[f,Y]+[f,-Y]=0$.
\end{step}

The cobordism group of immersions of oriented $5$-manifolds into $\R^8$ is isomorphic to
the $8^{\rm th}$ stable homotopy group of the Thom space $MSO(3)$ which is isomorphic to the $9^{\rm th}$ stable homotopy group of the suspension $\Sigma MSO(3)$. The map on the cobordism group induced by changing the orientation on the source coincides with the map on the stable homotopy group induced by the involution $\iota$ of $\Sigma MSO(3)$ which is reflection in $MSO(3)$. It follows from standard properties of suspensions that if $Z$ is any topological space and if $\img(\pi_n(Z))$ denotes the image of the suspension map $\pi_n(Z) \to \pi_{n+1}(\Sigma Z)$ then the map $\iota_\ast$, induced by the involution $\iota\colon\Sigma Z\to\Sigma Z$ which is reflection in $Z$, agrees with multiplication by $-1$ on $\img(\pi_n(Z))$. Thus there exists an oriented $6$-manifold $W$ with $\pa W=Y\sqcup (-Y)$ and an immersion $\tilde F\colon W\to\R^9_+$ which agrees with $\tilde f$ on both boundary components.

Let $\nu$ denote the normal bundle of $\tilde F$. Consider the manifold $\hat W=W\cup(Y\times I)$ obtained by identifying boundary components. Note that $\hat W$ is orientable. Moreover, we can extend $\nu$ to a bundle over $\hat W$ in a canonical way: use the identity transition function at both junctions. Let $\hat \nu$ denote this extension. Note that $T\hat W\oplus\hat\nu$ is a trivial bundle and thus $\hat \nu$ is a normal bundle for $\hat W.$ Consider a unit vector $v\in\R^8$ as a vector field along $\tilde F(W)$. Projecting $v$ to $\nu$ we get a section $n$ of the bundle $\nu$ over $W$. Note the section $n$ can be continued in a canonical way along $Y\times I$ and thus gives a section of $\hat \nu$. The zero set $\Sigma$ of $n$ is dual to $\w_3(\hat W)$. Moreover, along $\Sigma$, $v$ gives a vector field in the restriction $T_\Sigma \hat W$. A point where $v$ is tangent to $\Sigma\cap W$ corresponds to a cusp of $F$, where $F=\pi_v\circ\tilde F$ and $\pi_v$ is a projection parallel to $v$.

For generic $v$, $\pi_v\circ \tilde f$ does not have any cusps. Thus, for such $v$ the corresponding    section of $T_\Sigma\hat W$ is nowhere tangent to $\Sigma\cap Y\times I$. Furthermore the $\mod 2$ number of tangency points is equal to the $\mod 2$ self intersection number $[\Sigma]\bullet[\Sigma]$ of $\Sigma$ in $\hat W.$
By Poincar{\'e} duality this self-intersection number is
$$
[\Sigma]\bullet[\Sigma]={\w}_3^2[\hat W]=0,
$$
where the last equality holds since $\hat W$ is an orientable $6$-manifold and $\Omega_6=0$. We conclude that if $F=\pi_v\circ\tilde F$ then $F$ is a cobordism with an even number of cusps. Lemma \ref{l3} then implies that $[f,Y]+[f,-Y]=0$.

We are now in position to finish the proof. By Steps \ref{step1}--\,\ref{step3}, we have
$$
[f,Y]=-[r\circ f, -Y]=-[f,-Y]=-[f,Y].
$$
We conclude that $[f,Y]$ has order $2$ and hence that (ii) holds.
\end{proof}

\section{Fold bordism groups}\label{sec:bor-cob}
In this section we describe fold bordism groups. We first introduce some notation. Following Stong \cite{St2}, we let $C(n,k)$ denote the bordism group of maps of $n$-manifolds into $(n+k)$-manifolds. More precisely an element in $C(n,k)$ is an equivalence class of maps $f\colon V\to X$ where $V$ is a closed $n$-manifold and $X$ is a closed $(n+k)$-manifold, and two such maps $f_j\colon V_j\to X_j$, $j=0,1$, are equivalent if there exist cobordisms $W$ and $Y$ with $\pa W=V_0\sqcup V_1$ and $\pa Y = X_0\sqcup X_1$, respectively, and a map $F\colon W\to Y$ extending $f_j$, $j=0,1$. Addition in $C(n,k)$ is induced by disjoint union.
Analogously, we let $C_{\rm so}(n,k)$ denote the oriented bordism group. The definition is the same as for $C(n,k)$ except that all source- and target manifolds and all cobordisms are required to be oriented.

A slight generalization of the Thom-Pontryagin construction shows that
\begin{align*}
C(n,k) &\approx \N_{n+k}(\Omega^{\infty}MO(k+\infty)),\\
C_{\rm so}(n,k) &\approx \Omega_{n+k}(\Omega^{\infty}MSO(k+\infty)).
\end{align*}
Here $\infty$ should be understood as any sufficiently large integer, $\Omega^{j} X$ denotes the $j^{\rm th}$ loop space of $X$, $\N_j(X)$ denotes the $j^{\rm th}$ bordism group of $X$, and $\Omega_j(X)$ denotes the $j^{\rm th}$ oriented bordism group of $X$.

The fold cobordism groups studied in Section \ref{sec:2} are obtained form the corresponding bordism groups by imposing the additional constraint that the maps and the cobordisms have only $\Sigma^{1,0}$-singularities. Similarly we define the groups $C^{1,0}(n,k)$ and $C^{1,0}_{\rm so}(n,k)$ using the definition of the groups $C(n,k)$ and $C_{\rm so}(n,k)$ given above {\em with the additional constraint} that the maps ($f\colon V\to X$ above) and the bordisms ($F\colon W\to Y$ above) have only $\Sigma^{1,0}$-singularities.

The following theorem describes the bordism groups of fold maps in the largest codimension where transversality arguments do not imply an isomorphism of these groups with the corresponding (unrestricted) bordism groups computed by Stong.
\begin{theo}\label{t:last}
The bordism groups of fold maps are as follows.
\begin{itemize}
\item[{\rm (a)}] For any $k\ge 0$,
$$
C^{1,0}(2k+1, k) \approx C(2k+1,k),
$$
\item[{$\rm (b_0)$}]
$$
C^{1,0}_{\rm so}(1,0)\approx\Z_2
$$
\item[{$\rm (b_1)$}]
$$
C_{\rm so}^{1,0}(5,2)\approx C_{\rm so}(5,2)\approx \Z_2
$$
\item[{$\rm (b_\ast)$}] For any $m\ge 2$
$$
C^{1,0}_{\rm so}(4m+1,2m)\approx
C(4m+1,2m).
$$
\item[{\rm (c)}] For any $m\ge 1$
$$
C^{1,0}_{\rm so}(4m-1,2m-1)\approx C_{\rm so}(4m-1,2m-1)\oplus\Z_{3^u},
$$
where the power $u$ satisfies $0\le u\le t$, $t = \min\{j\mid \alpha_3(2m + j)\leq 3j\}$, see Theorem \ref{t1} (c).
\end{itemize}
\end{theo}

\begin{proof}
The forgetful maps
\begin{align*}
&r\colon C^{1,0}(2k+1,k)\to C(2k+1,k)\\
&r_{\rm so}\colon C^{1,0}_{\rm so}(2k+1,k)\to C_{\rm so}(2k+1,k),
\end{align*}
are obviously surjective. It is a consequence of Theorem \ref{l2} (a) and Lemma \ref{l3} that the map $r$ is also injective. This implies (a). Similarly ${\rm (b_\ast)}$ follows from Theorem \ref{l2} ${\rm (b_\ast)}$ in combination with Lemma \ref{l3}.

Consider ${\rm (b_0)}$. The number of cusps of a stable map from a closed oriented surface to an oriented surface is even and pairs can be removed by cobordism, see e.g. \cite{Ka}, and the group $C_{\rm so}(1,0)$ is trivial. An argument entirely analogous to the proof of Theorem \ref{t1} ${\rm (b_0)}$ now shows that ${\rm (b_0)}$ holds.

Consider ${\rm (b_1)}$. We construct a map $f\colon M\to\C P^4$ of an oriented $6$-manifold $M$ with the property
$$
\langle(w_3^2+w_4w_2)(\nu_f),[M]\rangle =1,
$$
where $\nu_f$ is the (virtual) normal bundle  of $f$. This implies that after small perturbation, making $f$ stable, it is a map with an odd number of cusps, see \eqref{eq:nofcusp}. With the existence of such a map established, the proof of ${\rm (b_1)}$ is a repetition the proof of Theorem 1 ${\rm (b_\ast)}$.

To construct such a map $f\colon M\to\C P^4$ it is sufficient to construct $F\colon\Sigma^L\C P^4\to MSO(L+2)$, where $\Sigma^L$ denotes the $L^{\rm th}$ suspension for some sufficiently large $L$, such that
\begin{equation}\label{e:F}
\langle F^\ast(Uw_3^2+Uw_4w_2),[\Sigma^L\C P^4]\rangle =1,
\end{equation}
where $U$ is the Thom class of $MSO(L+2)$. Indeed, applying the Thom-Pontryagin construction to $F$ gives an embedding $i\colon M\to\Sigma^L\C P^4$ of an orientable $6$-manifold $M$ such that
$$
(w_3^2+w_4w_2)(\nu_i)=(F|M)^\ast(w_3^2+w_4w_2),
$$
where the cohomology classes in the left hand side are the Steifel-Whitney classes of the universal $SO(L+2)$-bundle. Without loss of generality, we may assume that $i(M)\subset\R^L\times\C P^4\subset\Sigma^L\C P^4$. Defining $f=\pi\circ i$, where $\pi\colon\R^L\times\C P^4\to\C P^4$ is the natural projection, we obtain $f\colon M\to\C P^4$ with properties as desired. To see this, note that the virtual normal bundle  $\nu_f$ of $f$ belongs to the stable equivalence class of $\nu_i$ and hence
\begin{align*}
\langle(w_3^2+w_4w_2)(\nu_f),[M]\rangle&=\langle(w_3^2+w_4w_2)(\nu_i),[M]\rangle=
\langle (F|M)^\ast(w_3^2+w_4w_2),[M]\rangle\\
&=\langle F^\ast(Uw_3^2+Uw_4w_2),[\Sigma^L\C P^4]\rangle =1.
\end{align*}

To construct the map $F\colon \Sigma^L\C P^4\to MSO(L+2)$, we note that the dimension of $\Sigma^L\C P^4$ equals $L+8$. Thus we can replace $MSO(L+2)$ with the homotopically $(L+8)$-equivalent space $K=K(\Z,L+2)\times K(\Z,L+6)\times
K(\Z_2,L+7)$, where $K(G,n)$ denotes the Eilenberg-Maclane space with $\pi_n(K(G,n))\approx G$. A standard calculation shows that the class $Uw_3^2+Uw_4w_2\in H^{L+8}(MSO(L+2);\Z_2)$ is identified with the class
$Sq^4Sq^2l_{L+2}+Sq^2l_{L+6}+Sq^1l'_{L+7}\in H^{(L+8)}(K;\Z_2)$ under an $(L+8)$-equivalence, where $l_N\in
H^N(K(\Z,N);\Z_2)$ and $l'_N\in H^N(K(\Z_2,N);\Z_2)$ are the cohomological
fundamental classes, and where $Sq^n$ denotes the degree $n$ Steenrod operation. Define $F$ as the composition $j\circ\kappa$ of the natural inclusion $j\colon K(\Z,L+2)\to K$ and the map $\kappa\colon\Sigma^L\C P^4\to K(\Z,L+2)$ that corresponds to the generator $\Sigma^La$ of the group
$H^{L+2}(\Sigma^L\C P^4)\approx H^2(\C P^4)=\Z\langle a\rangle$ (i.e., the map $\kappa$ such that $\kappa^\ast l_{L+2}=\Sigma^L a$). The map $F$ then satisfies
\begin{align*}
\langle F^\ast(Uw_3^2+Uw_4w_2),[\Sigma^L\C P^4]\rangle&=
\langle F^\ast(Sq^4Sq^2l_{L+2}+Sq^2l_{L+6}+Sq^1l'_{L+7}),[\Sigma^L\C P^4]\rangle\\
&=\langle\kappa^\ast j^\ast(Sq^4Sq^2l_{L+2}+Sq^2l_{L+6}+Sq^1l'_{L+7}),[\Sigma^L\C P^4]\rangle\\
&=\langle\kappa^\ast j^\ast(Sq^4Sq^2l_{L+2}),[\Sigma^L\C P^4]\rangle\\
&=\langle Sq^4Sq^2\kappa^\ast j^\ast l_{L+2},[\Sigma^L\C P^4]\rangle =\langle Sq^4Sq^2\Sigma^La,[\Sigma^L\C
P^4]\rangle\\
&=\langle \Sigma^LSq^4Sq^2a,[\Sigma^L\C P^4]\rangle =\langle Sq^4Sq^2a,[\C P^4]\rangle\\
&=\langle Sq^4a^2,[\C P^4]\rangle=\langle a^4,[\C P^4]\rangle =1.
\end{align*}
Thus, \eqref{e:F} holds as desired.

Finally, consider (c). Let $G_m$ denote the greatest common divisor of the numbers $$\bigl\{{\overline p(f,Y)}_m[W]\bigr\},$$ where $W$ ranges over all closed oriented $4m$-manifolds, where $X$ ranges over all closed $(6m-1)$-manifolds, where $f$ ranges over all maps $f\colon W\to Y$, and where ${\overline p(f,Y)}_m$ denotes the $m^{\rm th}$ Pontryagin class of the virtual normal bundle $[f^\ast TY]-[TW]$ of $f$. It follows from Theorem \ref{l2} (c) that $G_m=3^u$ where $u$ satisfies the requirements in the statement. The argument used in the proof of Theorem \ref{t1} (c) allows us to define an homomorphism
$$
\krn\bigl(r_{\rm so}\colon C^{1,0}_{\rm so}(4m-1,2m-1)\to C_{\rm so}(4m-1,2m-1)\bigr)\to \Z_{G_m},
$$
by counting the algebraic number of cusps in a bordism from a given representative to the empty set. Again, using the map in Figure \ref{fig:cusp} and Lemma \ref{l3}, we conclude that this homomorphism is an isomorphism. This proves that there exists an exact sequence
$$
\begin{CD}
0 @>>> \Z_{3^u} @>>>C^{1,0}_{\rm so}(4m-1,2m-1) @>>>C_{\rm so}(4m-1,2m-1) @>>> 0.
\end{CD}
$$
This sequence splits since $C_{\rm so}(4m-1,2m-1)\approx\Omega_{6m-2}(\Omega^\infty MSO(2m-1+\infty)$ does not have odd torsion, see \cite{St3}, Chapter 4. Indeed, $\tilde H^{\ast+2m-1+N}(MSO(2m-1+N);\Z)\approx
H^\ast (BSO(2m-1+N);\Z)$, where $\tilde H$ denotes reduced cohomology, and the cohomology ring of the Grassmann-manifold does not have odd torsion, see \cite{T2}. Hence $\Omega_\ast(\Omega^N MSO(2m-1+N))$ does not have odd torsion either. This finishes the proof of (c).
\end{proof}

\section{Fold cobordism groups of general target manifolds}\label{sec:4}
In this section we describe fold cobordism groups of $(2k-1)$-manifolds into
an arbitrary closed $(3k-2)$-manifold $P$. More precisely, two fold maps
$f_j\colon M_j\to P$, $j=0,1$, of closed manifolds are {\em fold cobordant} if
there exists a cobordism $W$ and a fold map $F\colon W\to P\times[0,1]$
extending $f_j\times\{j\}$, $j=0,1$, on $\pa W$. In case $P$ is oriented we
define oriented fold cobordism analogously, requiring all manifolds and
cobordisms to be oriented. We denote the equivalence classes of fold cobordant
maps and oriented fold cobordant maps $\Sigma^{1,0}(2k-1,P)$ and $\Sigma_{\rm
  so}^{1,0}(2k-1,P)$, respectively. In contrast to the case $P=\R^{3k-2}$ the
set $\Sigma^{1,0}(2k-1,P)$ form a semi-group rather than
a group. (There is no natural geometric construction of an inverse operation.)
The following result shows that the semi-group structures on
$\Sigma^{1,0}(P,k)$ and $\Sigma^{1,0}_{\rm so}(P,k)$ in fact come from Abelian
group structures and moreover we describe the corresponding groups.

\begin{theo}\label{t:lastlast}
The cobordism semi-groups of fold maps into a manifold $P$ of dimension $(3k-2)$ satisfy the following.
\begin{itemize}
\item[{\rm (a)}] For any $k\ge 0$,
$$
\Sigma^{1,0}(2k+1, P) \approx \N_{2k+1}(P),
$$
\item[{$\rm (b_0)$}]
$$
\Sigma^{1,0}_{\rm so}(1,S^1)\approx\Z\oplus \Z_2
$$
\item[{$\rm (b_1)$}] If $P$ is a $7$-dimensional orientable manifold such that $H_1(P;\Z_2)\approx 0\approx H_2(P;\Z_2)$ then there exists an exact sequence of Abelian groups
$$
\begin{CD}
0 @>>> \Z_2 @>>> \Sigma_{\rm so}^{1,0}(5,P) @>>> \Omega_5(P) @>>> 0
\end{CD}
$$
\item[{$\rm (b_\ast)$}] For any $m\ge 2$ if $\dim(P)=6m+1$ and if $P$ is orientable then
$$
\Sigma^{1,0}_{\rm so}(2k+1, P) \approx \Omega_{2k+1}(P).
$$
\item[{\rm (c)}] For any $m\ge 1$, if $\dim(P)=6m-2$ and $P$ is orientable then there is an exact sequence of Abelian groups
$$
\begin{CD}
0 @>>> \Z_{3^v} @>>>\Sigma^{1,0}_{\rm so}(4m-1,P) @>>>\Omega_{4m-1}(P) @>>> 0,
\end{CD}
$$
where the power $v$ satisfies $0\le v\le t$, $t = \min\{j\mid \alpha_3(2m + j)\leq 3j\}$, see Theorem \ref{t1} (c).
\end{itemize}
\end{theo}

\begin{proof}
The proofs of (a) and ${\rm (b_\ast)}$ are completely analogous to the proofs of the corresponding parts of Theorems \ref{t1} and \ref{t:last}: Theorem \ref{l2} implies that any cobordism between two fold maps can be changed into a cobordism with an even number of cusps and Lemma \ref{l3} then implies we can cancel the cusps in pairs giving a fold cobordism.

The proof of ${\rm (b_0)}$ is a straightforward modification of the corresponding part of the proof of Theorem \ref{t1}, the $\Z$-summand comes from the mapping degree.

Consider ${\rm (b_1)}$. We first show that if
\begin{equation}\label{e:triv}
\Hom(\Omega_6(P),\Z_2)=0
\end{equation}
then there is an exact sequence as claimed. The condition \eqref{e:triv} implies that any stable map of a closed oriented $6$-manifold into $P\times I$ has an even number of cusps: the argument in the proof of Theorem \ref{l2} ${\rm (b_{01})}$ shows that the $\mod 2$ number of cusps of such maps gives a homomorphism $\Omega_6(P)\to\Z_2$. Therefore, as in proof of Theorem \ref{t1} ${\rm (b_1)}$, we can define a map $\Gamma$ on the preimage in $\Sigma^{1,0}_{\rm so}(5,P)$ of the neutral element of $\Omega_5(P)$ under the forgetful morphism into $\Z_2$ by counting the $\mod 2$ number of cusps in any cobordism of a fold map representing a given element. Moreover, Lemma \ref{l3} implies that this map is an isomorphism. This establishes the exact sequence.

We next claim that the homological condition on $P$ is equivalent to \eqref{e:triv}. We have
\begin{equation}\label{e:bord}
\Omega_6(P)\approx \oplus_j H_j(P;\Omega_{6-j})= H_1(P;\Z_2)\oplus H_2(P;\Z)\oplus H_6(P;\Z),
\end{equation}
modulo odd torsion, see e.g. \cite{CF}. In particular, \eqref{e:triv} holds provided any homomorphism from the right hand side of \eqref{e:bord} into $\Z_2$ is trivial.

Write $H_1(P;\Z)=F_1\oplus T_1^{\rm o}\oplus T_1^{\rm e}$, where $F_1$ is the free part which is a direct sum of infinite cyclic groups and where $T_1^{\rm o}$ and $T_1^{\rm e}$ are direct sums of finite cyclic groups of odd respectively even order. The universal coefficient theorem then implies that $H_1(P;\Z)\approx F_1\otimes\Z_2\oplus T_1^{\rm e}\otimes\Z_2$ and that $H^1(P;\Z)\approx F_1$.
Thus,
$$
\Hom(H_1(P;\Z_2),\Z_2)=0
$$
is equivalent to $F_1=0=T_1^{\rm e}$. Poincar{\'e} duality implies $H_6(P;\Z)\approx H^1(P;\Z)\approx F_1$. Hence $\Hom(H_6(P;\Z),\Z_2)=0$ if and only if  $F_1=0$. Finally, if $F_1=0=T_1^{\rm e}$ then $\Ext(H_1(P,\Z),\Z_2)=0$ and the universal coefficient theorem gives
$$
H^2(P;\Z_2)\approx \Hom(H_2(P;\Z),\Z_2)
$$
and the claim follows. To finish the proof of ${\rm (b_1)}$, we recall the general algebraic fact that if
$$
\begin{CD}
0 @>>> A @>>> B @>>> C @>>>  0
\end{CD}
$$
is an exact sequence of commutative semi-groups where $A$ and $C$ are Abelian groups and where the neutral element of $A$ maps to the neutral element of $B$ then $B$ is an Abelian group as well.

Finally, the proof of (c) uses the algebraic fact in combination with a
 repetition of the proof of Theorem \ref{t:last} (c) with the following
 alteration. The number $G_m=3^u$ in the proof of Theorem \ref{t:last} (c)
 should be replaced in the present case
 by the number $G_m=3^v$ which is the greatest common divisor of the numbers in the set $\bigl\{{\overline p(f,P\times I)}_m[W]\bigr\}$, where $W$ ranges over all oriented $4m$-manifolds and $f$ over all maps.
\end{proof}

\end{document}